\documentclass[a4paper,11pt]{article}
\pdfoutput=1
\usepackage{jheppub}
\usepackage{booktabs}
\usepackage{cases}

\newcommand{\ii}{\text{i}}
\newcommand{\dd}[1]{\text{d}#1 \, }
\newcommand{\e}[1]{\text{e}^{#1}}

\title{Aspects of Hecke Symmetry I: Ramanujan Identities and Inversion Formulas}

\author{Madhusudhan Raman}
\affiliation[a]{Institute of Mathematical Sciences,
Homi Bhabha National Institute (HBNI)\\
IV Cross Road, C.~I.~T.~Campus,
Taramani, Chennai 600 113,
Tamil Nadu, India}

\affiliation[b]{Department of Theoretical Physics,
Tata Institute of Fundamental Research\\
Homi Bhabha Road, Navy Nagar, Colaba, Mumbai 400 005,
Maharashtra, India}

\emailAdd{madhur@imsc.res.in}

\abstract{We study two aspects of Hecke symmetry in this note: first, we conjecture a generalization of the Ramanujan identities to the case of automorphic forms of Hecke groups; second, we conjecture a generalization of an inversion formula from the theory of elliptic functions.}

\begin{document} 
\maketitle
\flushbottom

\section{Introduction}
The goal of this paper is make a case for two conjectures: the first is a generalization of the Ramanujan identities to the case of automorphic forms; the second is a generalization of an inversion formula from the theory of elliptic functions. Both these conjectures will use Hecke symmetry in a crucial way.

Consider the following series:
\begin{align}
	\label{eq:RamanujanP}
	P(q) &= 1-24\sum_{k=1}^{\infty} \frac{k \, q^k}{1-q^k} \ , \\
	\label{eq:RamanujanQ}
	Q(q) &= 1+240\sum_{k=1}^{\infty} \frac{k^3 \, q^k}{1-q^k} \ , \\
	\label{eq:RamanujanR}
	R(q) &= 1-504\sum_{k=1}^{\infty} \frac{k^5 \, q^k}{1-q^k} \ .
\end{align}
Here, $q$ is the elliptic nome, defined as
\begin{equation}
q = \e{2\pi\ii \, \tau} \ ,
\end{equation}
where $\tau$ is a complex number valued in the upper-half plane $\mathbb{H}$. Expansions in the limit $ q \rightarrow 0 $ (equivalently, $ \tau \rightarrow \ii\infty $) of eqs.~\eqref{eq:RamanujanP} to \eqref{eq:RamanujanR} are Fourier expansions of Eisenstein series associated to SL$ (2,\mathbb{Z}) $ under the identifications $ \left(P, Q, R\right) \leftrightarrow \left(E_2, E_4, E_6\right) $.
The weight-$4$ ($E_4$) and weight-$6$ ($E_6$) Eisenstein series transform as
\begin{equation}
E_k \left( \frac{a \tau + b}{c \tau + d} \right) =( c \tau + d )^k \, E_k (\tau) \ .
\end{equation}
This is not true for $k=2$; instead, we find the anomalous transformation property
\begin{equation}
E_2 \left( \frac{a \tau + b}{c \tau + d} \right) =( c \tau + d )^2 \, E_2 (\tau) + \frac{6c}{\ii\pi} \, (c\tau+d) \ .
\end{equation}
The anomalous transformation property of $E_2$ is the reason it is referred to as a \emph{quasi-modular} form. The ring of quasi-modular forms is identified with $\mathbb{C}[E_2,E_4,E_6]$.

In 1916, Ramanujan \cite{Ramanujan1916} proved the following identities:
\begin{align}
	\label{eq:RamanujanDP}
	q \frac{\dd{P}}{\dd{q}} &= \frac{1}{12} \left( P^2 - Q\right) \ , \\
	\label{eq:RamanujanDQ}
	q \frac{\dd{Q}}{\dd{q}} &= \frac{1}{3} \left( PQ - R\right) \ , \\
	\label{eq:RamanujanDR}
	q \frac{\dd{R}}{\dd{q}} &= \frac{1}{2} \left( PR - Q^2\right) \ .
\end{align}
We have already seen that for the modular group SL$(2,\mathbb{Z})$, all quasi-modular forms may be expressed in terms of polynomials of $P$, $Q$, and $R$. Consequently, it is sufficient to specify eqs.~\eqref{eq:RamanujanDP} to \eqref{eq:RamanujanDR}, which may be used to determine the derivatives of all heavier forms---this is the statement that the ring of quasi-modular forms is closed under derivations. Now that we understand what we have set out to generalize, we will now turn our attention to the direction in which this generalization will proceed.

In an elegant note, Alperin \cite{Alperin} argued that
\begin{equation}
\text{PSL}(2,\mathbb{Z}) \cong C_2 * C_3 \ ,
\end{equation}
where $ C_k $ denotes a cyclic group of order $ k $ and $ * $ denotes the free product. This characterization admits a natural generalization to the case
\begin{equation}
\text{H}(m) = C_2 * C_m \ ,
\end{equation}
where $m\in\mathbb{Z}$ and $m \geq 3$. The groups H$(m)$ so constructed and indexed by $m$, which we will refer to as the \emph{height}, are called Hecke groups; they form an infinite family of discrete subgroups of SL$(2,\mathbb{R})$. With this embedding, our perspective has widened: we will treat SL$(2,\mathbb{Z})$ as the shortest of the Hecke groups, with height $m=3$. Above it stands an infinite tower of groups, each labelled by an integer height, and each with its own native automorphic forms---these are analogues of the Eisenstein series of SL$(2,\mathbb{Z})$. Our first goal in this paper will be to conjecture analogues of Ramanujan identities for automorphic forms of \emph{all} Hecke groups. 

Throughout this paper we will identify objects (functions, forms, etc.) by their $q$-expansions. This is perhaps closer to the spirit of Ramanujan's original work, although it is admittedly anachronistic.

\paragraph{Organization}
This paper is organized as follows. In Section \ref{sec:ModularAutomorphic} we identify Hecke groups as part of a larger family of discrete subgroups of SL$(2,\mathbb{R})$ called triangle groups. We then describe the manner in which their hauptmoduln (analogues of $j$-invariants) and automorphic forms (analogues of Eisenstein series) are constructed. In Section \ref{sec:RamanujanFormulas} we conjecture a generalization of the Ramanujan identities for all Hecke groups. This formula requires as input the height $m$ of the Hecke group and the weight $k$ of the automorphic form of interest. We use this to define an analogue of the Ramanujan-Serre derivative for Hecke groups. In Section \ref{sec:InversionFormula} we continue our exploration of Hecke symmetry, and use it to point out that an inversion formula from elliptic function theory admits a natural generalization. Finally, in Section \ref{sec:Discussion} we discuss the significance of this result for low-orders/low-orders resurgence in certain quantum mechanical systems.

\acknowledgments{It is our pleasure to thank Sujay Ashok, Shaun Cooper, Dileep Jatkar, Hossein Movasati, and P. N. Bala Subramanian for discussions and helpful correspondence. This work was done during a visit to the Tata Institute for Fundamental Research, Mumbai, and we acknowledge support from the Infosys Endowment for Research into the Quantum Structure of Spacetime. We are deeply grateful to Raju Bathija, Abhijit Gadde, Gautam Mandal, Shiraz Minwalla, Aniket Surve, and Sandip Trivedi for arranging our visit, and for their kindness and hospitality. Finally, we thank A. D. for encouragement and engaging conversations.}

\section{Modular and Automorphic Forms}
\label{sec:ModularAutomorphic}
We find it useful to understand Hecke groups as members of a wider class of discrete subgroups of PSL$(2,\mathbb{R})$ that are called \emph{triangle} groups. Our presentation in this section will closely follow the work of Doran et al.~\cite{Doran:2013npa}.

\subsection{Triangle and Hecke Groups}
For an integer triple $(m_1, m_2, m_3) \in \mathbb{Z}^+$, consider a triangle on a hyperbolic plane with angles $\pi/m_1$, $\pi/m_2$, and $\pi/m_3$. One can tesselate the hyperbolic (i.e.~upper-half) plane by reflecting such a triangle across its sides; the generators of these reflections form a discrete subgroup of PSL$(2,\mathbb{R})$ called a \emph{triangle group} $\Gamma_{\mathfrak{t}}$, specified by a \emph{type} $\mathfrak{t} = (m_1, m_2, m_3)$.

Hyperbolic geometry allows for triangles with \emph{cusps}: when the angle subtended by two of its sides is zero. In this paper we will restrict our attention to triangles with a single cusp, and by convention, we will realize this cusp by setting $m_3 = \infty$. This corresponds to triangle groups of type $(m_1, m_2, \infty)$. Finally, Hecke groups are a special case of this class with $m_1 = 2$ and $ m_2 \geq 3 $. Henceforth, we shall denote Hecke groups by the symbol H$(m) \cong \Gamma_{(2,m,\infty)}$. This uniquely defines Hecke groups up to conjugation by elements of PSL$(2,\mathbb{R})$. We will refer to $m$ as the \emph{height} of the Hecke group. The shortest Hecke group has a height $m=3$, and corresponds to the modular group. The tallest Hecke group has a height $m = \infty$, and corresponds to the congruence subgroup $\Gamma_0 (4)$ of the modular group. Our interest will be in all Hecke groups with finite heights, i.e.~all $3 \leq m < \infty$.

More concretely, we may define Hecke groups by their action: given a complex variable $\tau \in \mathbb{H}$, the Hecke group H$(m)$ is generated by the symbols $T$ and $S$ that act on it as
\begin{align}
	\label{eq:HeckeT}
	T \ &: \tau \longrightarrow \tau + 1 \ , \\
	\label{eq:HeckeS}
	S \ &: \tau \longrightarrow -\frac{1}{\lambda_m \tau} \ ,
\end{align}
where 
\begin{equation}
\label{eq:HeckeSpectrum}
\lambda_m = 4 \cos^2 \left( \frac{\pi}{m} \right) \ .
\end{equation}
The spectrum \eqref{eq:HeckeSpectrum} shows that there are precisely four cases where $\lambda_m \in \mathbb{Z}$, and we will refer to the corresponding groups as \emph{arithmetic} Hecke groups. Each of these arithmetic Hecke groups corresponds to a congruence subgroup $\Gamma_0(\lambda_m)$ of the modular group. Those H$(m)$ for which the corresponding $\lambda_m \notin \mathbb{Z}$ will be referred to as \emph{non-arithmetic} Hecke groups.

We observe that on setting $m=3$, the symbols in \eqref{eq:HeckeT} and \eqref{eq:HeckeS} will generate the modular group SL$(2,\mathbb{Z})$. For all integers $m > 3$, we obtain an infinite family of groups labelled by their heights. One of the goals of this paper is to investigate the extent to which constructions associated with the modular group admit generalizations to the Hecke case. In the following section, we review the construction of analogues of Klein's $j$-invariant for the Hecke groups.

\subsection{Hauptmoduln}
\label{sec:Hauptmoduln}
Given a Hecke group H$(m)$, there exists a unique weight-$0$ automorphic form---the hauptmodul---denoted $J_m$, which we will now construct. Let $q = \e{2\pi\ii\,\tau}$ as usual. The hauptmoduln corresponding to a Hecke group H$(m)$ has the following expansion about $\tau = \ii\infty$:
\begin{equation}
\label{eq:JAnsatz}
J_m = \frac{1}{\widetilde{q}} + \sum_{k=0}^{\infty} c_k \, \widetilde{q}^k \ ,
\end{equation}
where $\widetilde{q} = d_m \, q$; the quantity $d_m$ \cite[Theorem 1]{Doran:2013npa} is the value of the hauptmodul at the fixed point of the $S$-action, in much the same way that
\begin{equation}
j(\ii) = 1728 \ ,
\end{equation}
for the usual $j$-invariant of SL$(2,\mathbb{Z})$. While $ d_m $ is integer for the arithmetic subgroups, Wolfart \cite{Wolfart1983} has shown that $d_m$ corresponding to non-arithmetic Hecke groups are in fact transcendental. We will not require the actual values of $d_m$, and will concern ourselves solely with the $\widetilde{q}$-expansions.

The above expansion coefficients $c_k$ are determined \cite[Theorem 1]{Doran:2013npa} by a Schwarzian differential equation 
\begin{equation}
\label{eq:JSchwarzian}
-2 \, \dddot{J_m}\dot{J_m} + 3 \, \ddot{J_m}^2 = \dot{J_m}^4 \left( \frac{1-1/m^2}{J_m^2} + \frac{3/4}{\left(J_m-1\right)^2} + \frac{1/m^2 - 3/4}{J_m\left( J_m - 1\right)}\right) \ .
\end{equation}
Here, the derivatives are with respect to $2\pi\ii\,\tau$.

In the second part of this paper, it will be useful to normalize our hauptmoduln such that their corresponding Fourier series starts with $q^{-1}$. With this in mind, we define
\begin{equation}
\label{eq:SmallJBigJ}
j_m = d_m \, J_m \ .
\end{equation}

\subsection{Automorphic Forms}
\label{sec:AutomorphicForms}
In this section we construct automorphic forms of Hecke groups using the hauptmoduln we constructed in the previous section. For $k$ a positive, even integer such that $k \geq 4$ define an \emph{automorphic form} of weight $k$ following \cite[Theorem 2]{Doran:2013npa} as
\begin{equation}
\label{eq:WeightKAutomorphicForm}
f_{k} = (-1)^\frac{k}{2} \, \dot{J_m}^\frac{k}{2}  \, J_m^{\lceil \frac{k}{2m} \rceil - \frac{k}{2}} \, \left( J_m - 1 \right)^{\lceil \frac{k}{4} \rceil - \frac{k}{2}} \ .
\end{equation}
We will work with expressions such as these at the level of $\widetilde{q}$-series, i.e.~the above expression is taken to define the $\widetilde{q}$-series of an automorphic form of weight $k$ with respect to the group H$(m)$. In what follows, the appearance of an automorphic form $f_{k}$ should understood as corresponding to a particular Hecke group H$(m)$.

Automorphic forms of definite weight $k$ form a vector space, which we shall denote $\mathfrak{m}_k$. Define
\begin{equation}
\label{eq:DeltaK}
\delta_k = \frac{k}{2} - \left\lceil \frac{k}{4} \right\rceil - \left\lceil \frac{k}{2m} \right\rceil \ .
\end{equation}
Following \cite[Theorem 2]{Doran:2013npa}, a basis for the vector space $\mathfrak{m}_k$ of automorphic forms of weight-$k$ is given by the set
\begin{equation}
\left\lbrace f_{k,\ell} \cong f_k J_m^\ell \right\rbrace \quad \text{for} \quad 0 \leq \ell \leq \delta_k \ ,
\end{equation}
along with the convention that $f_k \equiv f_{k,0}$. Thus
\begin{equation}
\label{eq:DimensionMk}
\text{dim} \, \mathfrak{m}_k = \delta_k + 1 \ .
\end{equation}

We can also define quasi-automorphic forms for Hecke groups. Let $L= \text{lcm} \, (2,m)$, and define the automorphic discriminant $\Delta_m = f_{2L}$. We define a weight-$2$ quasi-automorphic form to be
\begin{equation}
f_2 = \frac{1}{2\pi\ii} \frac{\dd}{\dd \tau} \log \Delta_m \ ,
\end{equation}
Throughout this paper, we shall consider weight-$2$ quasi-automorphic forms $f_2$ to be normalized such that
\begin{equation}
\lim_{\tau\rightarrow\ii\infty} f_2 = 1 \ .
\end{equation}
It is easily verified that for $ m=3 $, the above prescriptions correctly reproduces the $q$-expansions \eqref{eq:RamanujanP}, \eqref{eq:RamanujanQ}, and \eqref{eq:RamanujanR} of the Eisenstein series $ E_2 $, $E_4$, and $E_6$ respectively. Further, for $ m > 3 $, the $ f_{k>6} $ are generically cusp-like forms, i.e.
\begin{equation}
\lim_{\tau\rightarrow\ii\infty} f_{k>6} = \mathcal{O}\left(q^{\delta_k}\right) \ ,
\end{equation}
and so they are not Eisenstein series, and should be thought of as an alternate choice of basis for the space of automorphic forms. Finally, for the arithmetic Hecke groups the reader might find it interesting to note that these weight-$2$ quasi-automorphic forms $f_2$ may be expressed in terms of the weight-$2$ quasi-modular Eisenstein series $E_2$ of SL$(2,\mathbb{Z})$ and certain quotients of Dedekind $\eta$-functions; see \cite[Appendix A]{Ashok:2016oyh} for explicit formulas. We now turn to the question of generalizing the Ramanujan identities to all Hecke groups.

\section{Ramanujan Identities}
\label{sec:RamanujanFormulas}
The automorphic forms corresponding to the Hecke groups, in particular their close structural similarity to modular forms corresponding to PSL$ (2,\mathbb{Z}) $, lend themselves to the following conjecture for generalizations of the Ramanujan identity.

\subsection*{Conjecture 3.1}
	The quasi-automorphic form $ f_2 $ corresponding to the Hecke group H$ (m) $  for all heights $3\leq m < \infty $ satisfies the following identity:
	\begin{equation}
	\label{eq:MainConjectureF2}
	f_{2}' =   \frac{1}{4} \left(\frac{m-2}{m}\right) \left(f_2^2 - f_4\right) \ ,
	\end{equation}

\subsection*{Conjecture 3.2}
	The automorphic forms $ f_k $ corresponding to the Hecke group H$ (m) $  for all heights $3\leq m < \infty $ satisfy the following identities for all weights $ 2 \leq k \leq m$:
	\begin{equation}
	\label{eq:MainConjecture}
	f_{2k}' =   \alpha \, f_2 f_{2k} + \beta \, f_{2k+2,\ell} + \gamma \, f_{2k+2} \ ,
	\end{equation}
	for the constants
	\begin{align}
		\begin{split}
			\alpha &= \frac{k}{2} \left( \frac{m-2}{m} \right) \ , \\
			\beta  &= \left( \left\lceil \frac{k+2}{2} \right\rceil - \left\lceil \frac{k+1}{2} \right\rceil \right) \left( \frac{2k-3m}{2m}\right) \ , \\
			\gamma &= (-1)^{k+1} \left( \frac{m-k}{m} \right) \ ,
		\end{split}
	\end{align}
	and the index $\ell$ defined such that
	\begin{equation}
	\ell = \begin{cases}
	0 & \text{if} \ k = m \ \text{and} \ m \in (2\mathbb{Z}+1) \ ,\\
	1 & \text{otherwise.}
	\end{cases}
	\end{equation}

These conjectures together form analogues of the Ramanujan identites for the case of all Hecke groups. This formula has been put to extensive tests in the parameter space $(m,k)$ of heights and weights, always yielding perfect agreement.

\paragraph{Ramanujan-Serre Derivatives}
Our relations above may be read as an analogue of the Ramanujan-Serre derivative  \cite{Zagier} for the case of Hecke groups. To define a derivation that sends forms in $\mathfrak{m}_{2k}$ to forms in $\mathfrak{m}_{2k+2}$, for $k>1$ we define a derivation
\begin{equation}
\text{D}: \mathfrak{m}_{2k} \longrightarrow \mathfrak{m}_{2k+2} \ ,
\end{equation}
whose action is effected by the following differential:
\begin{equation}
\text{D} = \frac{1}{2\pi\ii} \frac{\dd{}}{\dd{\tau}} - \frac{k}{2}\left( \frac{m-2}{m} \right)f_2 \ .
\end{equation}

\section{A Universal Inversion Formula}
\label{sec:InversionFormula}
There has been a sustained interest in inversion formulas in elliptic function theory starting with Ramanujan \cite{RamanujanNotebooks} and (independently) Fricke \cite{Fricke1916}, to Venkatachaliengar \cite{Venkatachaliengar}, Berndt et al.~\cite{BerndtBhargavaGarvan,BerndtChanLiaw}, Chan \cite{Chan98}, and more recently Cooper \cite{Cooper}. Our goal in this section is to point to a generalization of the work of these authors on inversion formulas for elliptic functions to alternative bases. We use the unifying perspective of Hecke groups and their hauptmoduln to conjecture this universal inversion formula. This will require the introduction of an infinite family of elliptic functions with \emph{rational} signatures; allowing the signatures to be rational greatly expands the class of alternative bases originally introduced by Ramanujan. As a consistency check, we reproduce (as special cases of our more general formula) four of Ramanujan's inversion formulas for elliptic functions to alternative bases.

\paragraph{Arithmetic Inversion Formula} We begin by reviewing earlier work. Following Cooper \cite{Cooper}, and for an index $r \in \lbrace 1, 2, 3, 4 \rbrace$, we define the function $x_r$ (for the case $r=1$) in terms of the Eisenstein series $Q$ and $R$ introduced in \eqref{eq:RamanujanQ} and \eqref{eq:RamanujanR} respectively:
\begin{equation}
\label{eq:CooperX1}
x_r(1-x_r) = \frac{Q^3 - R^2}{4\,Q^3} \ ,
\end{equation}
and for the cases $r = 2, 3,$ and $4$ by the following expression:
\begin{equation}
\label{eq:CooperX234}
\left( \frac{x_r}{1-x_r} \right)^{r-1} = r^6 q^{r-1} \prod_{k=1}^{\infty} \frac{\left( 1-q^{k r}\right)^{24}}{\left( 1-q^{k}\right)^{24}} \ .
\end{equation}
Here, the branches are determined by requiring that $x_r$ satisfy the following fall-offs near the weak-coupling ($\tau = \ii\infty$) point:
\begin{equation}
x_r = \begin{cases}
432 \, q + \mathcal{O}(q^2)  \ &\text{if} \ r=1 \ , \\
r^{\frac{6}{r-1}} \, q + \mathcal{O}(q^2) \ &\text{if} \ r=2, 3, \, \text{and} \, 4 \ .
\end{cases}
\end{equation}
The infinite product in \eqref{eq:CooperX234} can be condensed in a nice way. Consider the Dedekind $\eta$-function, defined as
\begin{equation}
\eta (\tau) = q^{1/24} \prod_{k=1}^{\infty} (1-q^k) \ ,
\end{equation}
and write \eqref{eq:CooperX234} in terms of an $\eta$-quotient like so:
\begin{equation}
\label{eq:CooperXEta}
\left( \frac{x_r}{1-x_r} \right)^{r-1} = r^6 \left( \frac{\eta(r \tau)}{\eta(\tau)}\right)^{24} \ .
\end{equation}
That \eqref{eq:CooperXEta} is naturally associated to the arithmetic Hecke groups $\Gamma_0 (r)$ is clear from the form of the $\eta$ quotient. However, is easy to see that this formula doesn't work for $r=1$, i.e. the case of SL$(2,\mathbb{Z})$.

The piecemeal nature of the definitions \eqref{eq:CooperX1} and \eqref{eq:CooperXEta} is unsatisfactory. As a step towards our generalization, we note that $r$ is nothing but $\lambda_m$ from \eqref{eq:HeckeSpectrum}, i.e.
\begin{equation}
\label{eq:HeckeSpectrum2}
r = 4 \, \cos^2 \left( \frac{\pi}{m} \right) \ ,
\end{equation}
If we are to generalize this discussion to include non-arithmetic Hecke groups, we will need to use the height $m$ to label these functions instead, with the correspondence between $r$ and $m$ given by \eqref{eq:HeckeSpectrum2}. Before we do so, let us introduce the subject of our generalization: the arithmetic inversion formula.

Consider definition of $x_r$ as presented in \eqref{eq:CooperX1}, which expresses $x_r$ as a function of $q$. The inversion formula in elliptic function theory corresponds to the statement that one can ``invert'' this formula to give $q$ as a function of $x_r$ instead. This is done in the following manner: for $0<x<1$, define $y$ as
\begin{equation}
y = \pi \, \frac{{}_2 F_1\left( \frac{1}{2},\frac{1}{2},1;1-x\right)}{{}_2 F_1\left( \frac{1}{2},\frac{1}{2},1;x\right)} \ .
\end{equation}
We observe that the function $y$ is simply the ratio of elliptic $K$ functions, defined as
\begin{equation}
K(x) = \frac{\pi}{2} \, {}_2 F_1 \left( \frac{1}{2},\frac{1}{2},1;x\right) \ ,
\end{equation}
and in terms of which $y$ is
\begin{equation}
y = \pi \frac{K(1-x)}{K(x)} \ .
\end{equation}
The inversion formula then states that the inversion of \eqref{eq:CooperX1} is given by
\begin{equation}
q = \e{-y} \ .
\end{equation}
Ramanujan \cite{Ramanujan1914} generalized this inversion formula by constructing elliptic functions to alternative bases:
\begin{equation}
\label{eq:rEllipticK}
K_r (x) = {}_2 F_1\left( a_r,1-a_r,1;x\right) \ ,
\end{equation}
for $a_r$ defined as
\begin{align}
	\label{eq:arCases}
	a_r = \begin{cases}
		\frac{1}{6} & \text{for} \ r = 1 \ , \\
		\frac{1}{4} & \text{for} \ r = 2 \ , \\
		\frac{1}{3} & \text{for} \ r = 3 \ , \\
		\frac{1}{2} & \text{for} \ r = 4 \ .
	\end{cases}
\end{align}
The quantity $a_r^{-1}$ is called the \emph{modular signature}. Observe that the modular signatures corresponding to these \emph{theories} are all integers. For $x_r$ defined by \eqref{eq:CooperX1} or \eqref{eq:CooperX234} as the case may be, we have the arithmetic inversion formula
\begin{equation}
q = \exp \left( - \frac{2\pi}{\sqrt{r}} \, \frac{K_r(1-x)}{K_r(x)} \right) \ .
\end{equation}
We now generalize this story to all Hecke groups. In order to do so, we will need to provide a definition of $ x_m $ that unifies the arithmetic cases \emph{and} generalizes straightforwardly to the non-arithmetic cases. 

\subsection*{Proposition 4.1}
	For all Hecke groups H$(m)$, define the function $x_m$ as
	\begin{equation}
	\label{eq:OurXAll}
	\frac{x_m}{x_m-1} = \frac{\sqrt{j_m - d_m}-\sqrt{j_m}}{\sqrt{j_m - d_m}+\sqrt{j_m}} \ .
	\end{equation}
	where $ j_m $ and $ d_m $ are introduced in Section \ref{sec:Hauptmoduln}. The branches are determined by the following fall-off condition near the weak-coupling point:
	\begin{equation}
	x_m = \frac{1}{4} \, \widetilde{q} + \mathcal{O}(\widetilde{q}^2) \ .
	\end{equation} 

It is straightforward to check that \eqref{eq:OurXAll} reproduces both \eqref{eq:CooperX1} and \eqref{eq:CooperX234} (see Ashok et al.~\cite{Ashok:2016oyh} for an explicit demonstration) in addition to admitting a generalization to arbitrary height $m \geq 3$. An explicit determination of the $ \widetilde{q} $-expansions confirms that this generalization satisfies all the required properties.

Our conjecture pertains to an inversion that gives $ q $ in terms of $ x_m $. We introduce elliptic functions with \emph{rational} signatures by generalizing \eqref{eq:rEllipticK} in the following manner:
\begin{equation}
K_m (x) = {}_2 F_1\left( a_m,1-a_m,1;x\right) \ ,
\end{equation}
where $a_m$ is defined as
\begin{equation}
a_m = \frac{m-2}{2m} \ .
\end{equation}
For the arithmetic cases $m \in \lbrace 3,4,6,\infty\rbrace$, $a_m$ reduces to the corresponding $a_r$ as defined in \eqref{eq:arCases}. Observe that for $m \notin \lbrace 3,4,6,\infty\rbrace$, we get modular signatures $a_m^{-1}$ that are \emph{rational}.

\subsection*{Conjecture 4.2}
	The formula \eqref{eq:OurXAll} is inverted by
	\begin{equation}
	\label{eq:InversionAll}
	q = \exp \left( - \frac{\pi}{\cos \left( \frac{\pi}{m} \right)} \, \frac{K_m(1-x_m)}{K_m(x_m)} \right) \ .
	\end{equation}
Once again, this inversion formula has been tested extensively against many heights $m$, and we have found complete agreement in all cases.

\section{Discussion}
\label{sec:Discussion}
Our interest in the properties of automorphic forms of Hecke groups and inversion formulas is their natural appearance in the theory of low-orders/low-orders (or P/NP) resurgence---see Dunne and Unsal \cite{Dunne:2013ada,Dunne:2014bca,Dunne:2016qix,Dunne:2016jsr}, and Codesido et al.~\cite{Codesido:2016dld,Codesido:2017jwp} for more details. In particular, our interest is in the Chebyshev oscillators, first studied in the beautiful work of Basar et al.~\cite{Basar:2017hpr}, where a version of \eqref{eq:InversionAll} first appeared. The authors there studied the classical and quantum mechanics of oscillators governed by the potentials
\begin{equation}
V(x) = T_{\mu}^2 (x) \ ,
\end{equation}
where $T_\mu(\cos \theta) = \cos \, \mu\theta$ is a Chebyshev polynomial of the first kind, and $\mu$ takes values in $\frac{1}{2}\mathbb{Z}$.\footnote{We use the symbol $\mu$ to avoid confusion with the height $m$ of the Hecke group.} They found that the special cases $\mu = \lbrace 3/2,2,3,\infty \rbrace$ inherit a modular structure in one-to-one correspondence with Ramanujan's theories of elliptic functions to alternative bases. Inspired by this observation, the author and P. N. Bala Subramanian \cite{Unpublished} showed that the identification of a Hecke symmetry underlying the classical mechanics of each Chebyshev oscillator \emph{persists} even for non-arithmeric Hecke groups. That is, to an oscillator labelled by $\mu$, we associate a (not necessarily arithmetic!) Hecke group with height 
\begin{equation}
m = 2\mu \ .
\end{equation}
Using this perspective, we were able to show that many of the formulas for classical observables---the periods and frequencies, for example---could be expressed in a unified manner in terms of automorphic forms of H$(2\mu)$ just like the ones we have encountered in this paper.

Observe that for $\mu \notin \lbrace 3/2,2,3,\infty \rbrace$, the classical relation between energy $u$ and momentum $p$ takes the form
\begin{equation}
\frac{1}{2} \, p^2 = u - T_\mu^2(x) \ ,
\end{equation}
which is a curve of genus $>1$. The question of whether these higher-genus potentials exhibit P/NP resurgence remains an interesting open problem, and it is our hope that a better understanding of their underlying automorphic structure might shed light on this issue. Our investigations here may be viewed as a small step in that direction.

Whether these formulas deserve the attention of mathematicians is a more difficult question to answer; there are, however, some interesting directions to pursue. For example, Ramamani \cite{Ramamani1970Thesis, Ramamani1989} and Toh \cite{Toh2011} have considered Ramanujan-like systems of differential equations for Eisenstein series corresponding to $\Gamma_0(2)$, and it would be interesting to see how our own work is related to theirs.

\paragraph{Note Added} After this work was completed, we learned of the work of Edvardsson \cite{Edvardsson}, whose results are obtained as special cases of Conjecture 3.2. A proof of an analogue of Conjectures $3.1$ and $ 3.2 $ for Eisenstein series has been furnished by the author and collaborators, see Ashok et al.~\cite{Ashok:2018myo} for more details.

\appendix

\appendix

\bibliographystyle{JHEP}
\bibliography{Refs}
\end{document}